\documentclass{article}
\usepackage{fullpage}
\def\C{{\rm C \kern-.48em\vrule width.06em height.57em depth-.02em \kern.48em}}
\def\Z{{{\rm Z}\kern-.28em{\rm Z}}}
\def\N{{{\rm I}\kern-.16em{\rm N}}}
\def\R{{{\rm I}\kern-.16em{\rm R}}}

\def\diag{\mathop{\rm diag}\nolimits}

\def\qRa{{{Q}{R}{A}}}

\def\eqbd{\mathop{{:}{=}}}
\def\bdeq{\mathop{{=}{:}}}
\newcommand {\eproof}
      {\hskip 4cm \hfill
        {\ \vbox{\hrule\hbox{\vrule height1.3ex\hskip0.8ex\vrule}\hrule}}
        \newline \vskip 0.3cm}

\begin{document}
\title{Evaluation of Sylvester type determinants using 
block-triangularization}

\author{Olga Holtz \thanks{
         Institut f\"ur Mathematik, MA 4-5,
         Technische Universit\"at Berlin,
         Stra{\ss}e des 17.\ Juni 136,
         D-10623 Berlin,
         Germany. \hskip 1cm  E-mail: {\tt{holtz@math.TU-Berlin.DE\/}}.
Supported by Alexander von Humboldt Foundation and by the DFG research center
``Mathematics for key technologies'' in Berlin.
Address at the time of publication: University of California-Berkeley,
Department of Mathematics, 821 Evans Hall, Berkeley, CA 94720, USA.}}
\date{}
\maketitle

\begin{abstract}
It is shown that the values of Sylvester type determinants for
various orthogonal polynomials considered 
by Askey in~\cite{A} can be ascertained inductively using simple 
block-triangularization schemes. 
\end{abstract}

\section{Introduction}

Richard Askey in~\cite{A} shows two ways, one matrix-theoretic and 
another based on orthogonal polynomials, to evaluate determinants 
\[  D_{N+1} (x) \eqbd \left |
    \begin{array}{cccccccc}
  x & 1 & &                &     &    &      \\
  N & x & 2 &            &    &  0  &     \\
    & N - 1 & x   &   3      &     &    &     \\
    &   & \cdots    &      &    &    &    \\
    & 0  &     &           & 2   & x   & N   \\
    &   &     &            &  0    & 1   & x   \\
     \end{array}
   \right | , \]
which were first considered by Sylvester~\cite{S}. In addition, he obtains
several generalizations of Sylvester's determinants and explores their 
connection to orthogonal polynomials.

The purpose of this note is to show how the determinants from~\cite{A}
can be evaluated in yet another way, based on partial information about 
left or right eigenvectors of the corresponding matrices coupled with a 
simple similarity trick. In all cases except one, only one (the most 
obvious) eigenvector is used. The exceptional case is the Sylvester 
determinant itself, where two eigenvectors are readily available and
hence used to derive the result.

\section{Sylvester's determinant and two close variants}

Let us start with the Sylvester determinant. We want to prove that
\begin{equation}
D_{N+1}(x)=\prod_{j=0}^N (x+N-2j)  \label{sylv}
\end{equation}
(which is formulas~(2.3),~(2.4) from~\cite{A}).
Since the values of $D_1$ and $D_2$ agree with~(\ref{sylv}), 
it is enough to show that
\begin{equation}
D_{N+1}(x)=(x-N)(x+N) D_{N-1}(x)  \label{indsylv}
\end{equation} 
for $N\geq 1$. The Sylvester determinant $D_{N+1}$ is the 
characteristic polynomial of the matrix $-{\cal D}_{N+1}$, i.e.,  
$D_{N+1}(x)=\prod_{\lambda \in \sigma(D_{N+1})} (x+\lambda)$, where 
\[ {\cal D}_{N+1} \eqbd \left[ \begin{array}{ccccccc}
0 & 1 & 0 & 0 & \cdots &0  & 0 \\ 
N & 0 & 2 & 0 & \cdots & 0 & 0\\
0 & N-1 & 0 & 3 & \cdots & 0 & 0 \\
0 & 0 & N-2 & 0 &  \cdots & 0 & 0 \\
\vdots & \vdots & \vdots & \vdots & \ddots & \vdots & \vdots \\
0 & 0 & 0 & 0 &  \cdots & 0 & N \\
0 & 0 & 0 & 0 &  \cdots & 1 & 0 \end{array}
\right]. \]
Note that $(1,1,1,\ldots, 1)$ is a left eigenvector of ${\cal D}_{N+1}$
corresponding to eigenvalue $N$ and $(1,-1,1,\ldots, (-1)^{N})$ is its left 
eigenvector corresponding to eigenvalue $-N$. The similarity transformation 
${\cal D}_{N+1} \mapsto {\cal T}_{N+1} {\cal D}_{N+1} {\cal T}_{N+1}^{-1}$ where
\[ {\cal T}_{N+1}\eqbd \left[ \begin{array}{rrrrrc} 1 & 1 & 1 & 1 & \cdots & 1 \\  
1 & -1 & 1 & -1 & \cdots & (-1)^N \\ 
0 & 0 & 1 & 0 & \cdots & 0 \\
0 & 0 & 0 & 1 & \cdots & 0 \\
\vdots & \vdots & \vdots & \vdots & \ddots & \vdots \\
0 & 0 & 0 & 0 & \cdots & 1 \end{array} \right] \]
therefore reduces ${\cal D}_{N+1}$ to a block lower-triangular form 
\[ \left[ \begin{array}{cc} \diag(N,-N) & 0 \\ {*} & {\cal M}_{N-1}
\end{array} \right], \]  where 
\[  {\cal M}_{N-1} \eqbd \left[ \begin{array}{ccccccc} 
0 & 3+(-N+1) & 0 &  -N+1 & 0 & -N+1 & \cdots  \\
 N-2 & 0 & 4 & 0 & 0 & 0 &  \cdots  \\ 
0 & N-3 &  0 & 5 & 0 & 0 & \cdots  \\
 0 & 0 & N-4 & 0 & 6 & 0 &  \cdots  \\ 
  0 & 0 & 0 & N-5 & 0 & 7 &  \cdots  \\ 
0 & 0 & 0 & 0 & N-6 & 0 &  \cdots  \\ 
\vdots  & \vdots &\vdots  & \vdots & \vdots & \vdots & 
\ddots  \end{array} \right].\]
Now it remains to show that ${\cal M}_{N-1}$ is similar to ${\cal D}_{N-1}$.
But this is indeed so, since ${\cal S}_{N-1}^{-1} {\cal M}_{N-1} {\cal S}_{N-1}=
{\cal D}_{N-1}$
where
\[ {\cal S}_{N-1} \eqbd \left[ \begin{array}{rrrrrr} 
1 & 0 & -1 & 0 & 0 &   \cdots  \\
0 & 1 & 0 & -1 & 0 & \cdots  \\ 
0 & 0 & 1 & 0 & -1 &  \cdots  \\ 
0 & 0 & 0 & 1 & 0 &   \cdots  \\ 
0 & 0 & 0 & 0 & 1 & \cdots  \\  
\vdots & \vdots & \vdots  & \vdots & \vdots & \ddots 
\end{array}  \right]. \]
This proves~(\ref{indsylv}) and therefore~(\ref{sylv}).
\eproof

A related determinant
\[ B_{N+1}  (x) \eqbd  \left|
     \begin{array}{cccccccc}
x        & a         & 0         & 0   &  \cdots & 0    & 0         & 0\\
N(a-1)  & x-1       & 2a        & 0   &  \cdots & 0    & 0         & 0\\
0        & (N-1)(a-1)& x-2       & 3a  &  \cdots & 0    & 0         & 0 \\                 &           & \cdots    &     &         &      &           &   \\
0        & 0         & 0         & 0   &  \cdots & 2(a-1) & x-(N-1) & Na  \\
0        & 0         & 0         & 0   &  \cdots &   0   & a-1       & x-N 
\end{array}    \right |  \]
can be evaluated analogously and even more simply. We use the notation
\[ {\cal B}_{N+1}  \eqbd  \left[
     \begin{array}{cccccccc}
0        & a         & 0         & 0   &  \cdots & 0    & 0         & 0\\
N(a-1)  & -1       & 2a        & 0   &  \cdots & 0    & 0         & 0\\
0        & (N-1)(a-1)& -2       & 3a  &  \cdots & 0    & 0         & 0 \\                 &           & \cdots    &     &         &      &           &   \\
0        & 0         & 0         & 0   &  \cdots & 2(a-1) & -(N-1) & Na  \\
0        & 0         & 0         & 0   &  \cdots &   0   & a-1       & -N 
\end{array}    \right]   \]
for the matrix that satisfies $B_{N+1}(x)=\prod_{\lambda \in 
\sigma({\cal B}_{N+1})} (x+\lambda)$. 

To prove that
\begin{equation}
B_{N+1}(x)=\prod_{j=0}^N [x+(N-2j)a-N+j] \label{sylvb} 
\end{equation}
(formula~(2.8) from~\cite{A}), we will show that
\begin{equation}
B_{N+1}(x)=(x+Na-N)B_N(x-a). \label{indsylvb} 
\end{equation}
The vector $(1, 1, 1, \ldots, 1)$ is a left eigenvector of ${\cal B}_{N+1}$
corresponding to eigenvalue $Na-N$. The similarity transformation
${\cal B}_{N+1} \mapsto {\cal T}_{N+1} {\cal B}_{N+1} {\cal T}_{N+1}^{-1}$ where
\[ {\cal T}_{N+1}\eqbd \left[ \begin{array}{cccccc} 1 & 1 & 1 & 1 & \cdots & 1 \\  
0 & 1 & 0 & 0 & \cdots & 0 \\ 
0 & 0 & 1 & 0 & \cdots & 0 \\
0 & 0 & 0 & 1 & \cdots & 0 \\
\vdots & \vdots & \vdots & \vdots & \ddots & \vdots \\
0 & 0 & 0 & 0 & \cdots & 1 \end{array} \right] \]
reduces ${\cal B}_{N+1}$ to a block lower triangular form
\[ \left[ \begin{array}{cc} Na-N & 0 \\ {*} & {\cal M}_N
\end{array} \right], \]  where 
\[  {\cal M}_N \eqbd \left[ \begin{array}{cccccc} 
-1-N(a-1) & 2a-N(a-1) & -N(a-1) &  -N(a-1) &  -N(a-1) & \cdots  \\
 (N-1)(a-1) & -2 & 3a & 0 & 0 &   \cdots  \\ 
0 & (N-2)(a-1) &  -3 & 4a & 0 & \cdots  \\
  0 & 0 & (N-3)(a-1) & -4 & 5a &  \cdots  \\ 
0 & 0 & 0 & (N-4)(a-1) & -5 &   \cdots  \\ 
\vdots  & \vdots &\vdots  & \vdots & \vdots & 
\ddots  \end{array} \right].\]
Finally, ${\cal M}_N$ is similar to ${\cal B}_N-aI$:
\[ {\cal M}_N {\cal S}_N={\cal S}_N ({\cal B}_N-aI),  \]
where
\[ {\cal S}_N\eqbd \left[ \begin{array}{rrrrrr} 1 & -1 & 0 & 0 & \cdots & 0 \\  
0 & 1 & -1 & 0 & \cdots & 0 \\ 
0 & 0 & 1 & -1 & \cdots & 0 \\
0 & 0 & 0 & 1 & \cdots & 0 \\
\vdots & \vdots & \vdots & \vdots & \ddots & \vdots \\
0 & 0 & 0 & 0 & \cdots & 1 \end{array} \right]. \] 
This proves~(\ref{sylvb}). \eproof

As is shown in~\cite{A}, the formula~(\ref{sylvb}) for $B_{N+1}$ implies
the formula~(\ref{sylv}) for the Sylvester determinant $D_{N+1}$, since
\[ D_{N+1}(x)=\lim_{a\to \infty} {B_{N+1}(ax) \over a^{N+1} }.  \]
Also,~(\ref{sylvb}) gives a formula for another related determinant,
\[ A_{N+1}  (x) \eqbd
\left | 
     \begin{array}{cccccccc}
x   & 1    & 0         & 0   &  \cdots & 0    & 0         & 0\\
-N  & x-2  & 2         & 0   &  \cdots & 0    & 0         & 0\\
0   &  -(N-1)     & x-4   &  3   &         &      &   &   \\
    & \cdots     &    &     &         &      &           &   \\
0   & 0    & 0         & 0   &         & -2   & x-2(N-1)  & N  \\
0   & 0    & 0         & 0   &  \cdots &      & -1        & x-2N
\end{array} 
       \right | \label{eq2.5} \]
via the relation
\[ A_{N+1}(x)=\left. 2^{N+1} B_{N+1}\left({x\over 2}\right) 
\right|_{a={1\over 2}}. \]

\section{Determinants for Krawtchouk and dual Hahn polynomials}

The determinant for the Krawtchouk polynomial $K_{N+1}(x;p,N)$ is given by
{\small
\begin{eqnarray*} K_{N+1}(x;p,N) \eqbd  \hbox{\hskip 13.5cm}  \\ 
  \left |
    \begin{array}{ccccccc}
  -x + p N   &  p N  &  0   & 0  &  \cdots  &  0  & 0  \\
  &&&&&&\\
   (1-p)     &  -x + p N + (1 - 2p)   &  p(N-1)   &  0  &  \cdots  &  0  & 0   \\
&&&&&&\\
   0        &    2(1-p)   & -x + p N + 2(1-2p)   &  p(N-2) &  \cdots  &  0  & 0  \\
&&&&&&\\
   &             &  \cdots               &         &          &     &   \\
&&&&&&\\
   0        &     0       &   0                  &  0      & \cdots   & N(1-p)  & -x + p N + N(1-2p)
      \end{array}
      \right |.  \end{eqnarray*} }
With the usual notation $(a)_k$ for the shifted factorial 
(see~\cite[formula~(3.10)]{A}), we need to prove~\cite[(3.25)]{A}:  
\begin{equation}
 K_{N+1}(x;p,N)=(-x)_{N+1}. \label{kraw} 
\end{equation}
For that, it is enough to show that
\begin{equation}
 K_{N+1}(x;p,N)=(-x)K_N(x;p,N-1). \label{indkraw} 
\end{equation}
As usual, we need to find an `obvious' eigenvector of the matrix
\[ {\cal K}_{N+1}\eqbd    
\left[  \begin{array}{ccccccc}
   p N   &  p N  &  0   & 0  &  \cdots  &  0  & 0  \\ 
&&&&&&\\
   (1-p)     &   p N + (1 - 2p)   &  p(N-1)   &  0  &  \cdots  &  0  & 0   \\
&&&&&&\\
   0        &    2(1-p)   &  p N + 2(1-2p)   &  p(N-2) &  \cdots  &  0  & 0  \\
&&&&&&\\
   &             &  \cdots               &         &          &     &   \\
&&&&&&\\
   0        &     0       &   0                  &  0      & \cdots   & N(1-p)  &  p N + N(1-2p)  \end{array}
      \right] \]
corresponding to eigenvalue $0$. In this case, it happens to be 
the right eigenvector $(1, -1, 1, -1, \ldots, (-1)^N)^T$. We then use
the transformation
\begin{equation}
 {\cal T}_{N+1}\eqbd \left[ \begin{array}{rccccc} 1 \quad  & 0 & 0 & 0 & \cdots & 0 \\  
-1\quad & 1 & 0 & 0 & \cdots & 0 \\ 
1\quad & 0 & 1 & 0 & \cdots & 0 \\
-1\quad & 0 & 0 & 1 & \cdots & 0 \\
\vdots\quad  & \vdots & \vdots & \vdots & \ddots & \vdots \\
(-1)^N & 0 & 0 & 0 & \cdots & 1 \end{array} \right] \label{t}
  \end{equation}
to transform ${\cal K}_{N+1}$ into a block upper-triangular form
\[ {\cal T}_{N+1}^{-1}{\cal K}_{N+1} {\cal T}_{N+1}  =
\left[ \begin{array}{cc} 0 & * \\ 0 & {\cal M}_N
\end{array} \right], \]  where 
 \[ 
{\cal M}_N\eqbd  \left[ \begin{array}{cccccc}
pN+PN+(1-2p) & p(N-1) & 0 & \cdots & 0 & 0 \\
-pN + 2(1-p) & pN+2(1-2p) & p(N-2) & \cdots & 0 & 0 \\
pN & 0 & 3(1-p)  & \cdots & 0 & 0 \\
-pN & 0 & 0  & \cdots & 0 & 0 \\
\vdots & \vdots & \vdots & \ddots & \vdots & \vdots \\
(-1)^N pN & 0 & 0 & \cdots & N(1-p) & pN+N(1-2p) \end{array}  
\right]. \]
Now,  ${\cal M}_N$ is similar to the matrix ${\cal K}_N+I$,
viz. 
\[{\cal S}_N{\cal M}_N =({\cal K}_N+I){\cal S}_N \]
via the transformation
\begin{equation}
 {\cal S}_N\eqbd \left[ \begin{array}{cccccc} 1 & 0 & 0 & 0 & \cdots & 0 \\  
1 & 1 & 0 & 0 & \cdots & 0 \\ 
0 & 1 & 1 & 0 & \cdots & 0 \\
0 & 0 & 1 & 1 & \cdots & 0 \\
\vdots & \vdots & \vdots & \vdots & \ddots & \vdots \\
0 & 0 & 0 & 0 & \cdots & 1 \end{array} \right]. \label{s}
\end{equation}
This proves~(\ref{indkraw}) and hence~(\ref{kraw}). \eproof

The determinant for dual Hahn polynomials,
{\small \begin{eqnarray*}
 R_{N+1}(\lambda(x);\gamma , \delta,N)  \eqbd \hbox{\hskip 13cm} \\
 \left | 
  \begin{array} {ccccc}
    \lambda (x) + N(\gamma + 1) & N(\gamma + 1) & 0 & 0  &     \\
                                &               &   &    &    \\
    (N+\delta )   &
    \begin{array} {r}
      \lambda (x) + N(\gamma + 3) \\
          - (\gamma - \delta + 2)
    \end{array}
        & (N - 1)(\gamma + 2)  &   0  &  \cdots    \\
         &&&& \\
     0  &  2(N+\delta - 1)  &
     \begin{array}{r}
     \lambda (x) + N(\gamma + 5)\\
         -2(\gamma - \delta + 4)
      \end{array}
      & (N - 2)(\gamma + 3)     &\\
         &&&& \\
     &    &  \cdots     &    &    \\
    &&&&\\
    0 & 0 & \cdots &  N(\delta + 1) &
    \left \{ \begin{array}{r}
    \lambda (x) + N(2N + \gamma + 1) \\
     - N (2N + \gamma - \delta ) 
    \end{array} \right \}
    \end{array}   
      \right |, 
 \end{eqnarray*}}
where \[ \lambda(x)\eqbd -x(x+\gamma+\delta+1),\]
appears in~\cite{A} as the first example beyond Krawtchouk polynomials 
in terms of simplicity of recurrence coefficients. This determinant 
seems much harder to evaluate, but this in fact requires 
just an additional shift of parameters.  The formula for $R_{N+1}$ 
is~\cite[(4.5)]{A}:
\begin{equation} R_{N+1}(\lambda(x);\gamma , \delta,N)=(-x)_N 
(x+\gamma+\delta+1)_{N+1}.  \label{hahn}
\end{equation} 
It can be proved from the relation
\begin{equation}
 R_{N+1}(\lambda(x);\gamma , \delta,N)=\lambda(x)  
R_N(\lambda(x)+(\gamma+\delta+2);\gamma+1 , \delta+1,N-1). \label{indhahn}  
\end{equation}
The latter can be checked exactly as above, i.e., by applying the transformation 
${\cal T}_{N+1}$ of the form~(\ref{t}) to the matrix
\begin{eqnarray*}
 {\cal R}_{N+1}(\gamma , \delta)  \eqbd  \hbox{\hskip 13cm} \\
 \left [
  \begin{array} {ccccc}
     N(\gamma + 1) & N(\gamma + 1) & 0 & 0 &     \\
                                &  &   &   &    \\
    (N+\delta )   &
    \begin{array} {r}
       N(\gamma + 3) \\
          - (\gamma - \delta + 2)
    \end{array}
        & (N - 1)(\gamma + 2)   &  0 & \cdots    \\
         &&&& \\
     0  &  2(N+\delta - 1)  &
     \begin{array}{r}
     N(\gamma + 5)\\
         -2(\gamma - \delta + 4)
      \end{array}
     & (N - 2)(\gamma + 3)     &\\
         &&&& \\
    
    &    &  \cdots     &    &    \\
    &&&&\\
     
    0 & 0 & \cdots &N(\delta + 1) &
    \left \{ \begin{array}{r}
    N(2N + \gamma + 1) \\
     - N (2N + \gamma - \delta ) 
    \end{array} \right \}
    \end{array}   
      \right ], 
 \end{eqnarray*}
obtaining a block upper-triangular matrix 
\[ {\cal T}_{N+1}^{-1} {\cal R}_{N+1}(\gamma , \delta) {\cal T}_{N+1}  =
 \left[ \begin{array}{cc} 0 & * \\ 0 & {\cal M}_N
\end{array} \right], \]  where 

\begin{eqnarray*} 
 {\cal M}_N(\gamma, \delta) \eqbd \hbox{\hskip 13cm} \\ 
\left[ \begin{array}{ccccc}
{\begin{array}{r} N(\gamma+1)+N(\gamma+3) \\
-(\gamma-\delta-2)\end{array} }  & (N-1)(\gamma+2) & 0 & 0 & \cdots  \\
&&&& \\
{\begin{array}{r} -N(\gamma+1) \\ +2(N+\delta-1)
\end{array}} & {\begin{array}{r} N(\gamma+5) \\ 
-2(\gamma-\delta+4)\end{array}}  & (N-2)(\gamma+3)  
& 0 & \cdots  \\
&&&& \\
N(\gamma+1) & 3(N+\delta-2) & {\begin{array}{r} N(\gamma+7)  \\
-3(\gamma-\delta+6) \end{array}} & (N-3)(\gamma+4) &
\cdots  \\
&&&& \\
-N(\gamma+1) & 0 & 4(N+\delta+3) & {\begin{array}{r} N(\gamma+9) \\ 
-4(\gamma-\delta+8) \end{array}}  & \cdots  \\ 
 \vdots & \vdots & \vdots & \vdots & \ddots  \end{array} \right], 
\end{eqnarray*}
and then transforming ${\cal M}_N(\gamma, \delta)$ into 
${\cal R}_N (\gamma+1 , \delta+1)+(\gamma+\delta+2)I$ by
\[ {\cal M}_N(\gamma, \delta) \mapsto {\cal S}_N{\cal M}_N(\gamma, \delta)
{\cal S}_N^{-1},\]
where ${\cal S}_N$ has the form~(\ref{s}). This proves~(\ref{indhahn}) and~(\ref{hahn}).
\eproof

\section{Determinants for Hahn, Racah, and $q$-Racah polynomials}

The last examples appearing in~\cite{A} are three types of orthogonal 
polynomials, whose recurrence coefficients are fractional rather than 
polynomial. Recall that the recurrence coefficients appear as entries 
on the main diagonal and the first sub- and super- diagonals of the 
tridiagonal determinant that gives the value of the corresponding 
polynomial. The alternating sums of these coefficients are still 
zero, just like in the two examples we just considered in section~3. 
In other words, the determinants we now consider all have the form
\begin{equation}  \left |    \begin{array}{ccccccc}
\lambda (x) - b_0        & a_0         & 0         &   \cdots & 0    & 0         & 0\\
c_1  & \lambda (x) - b_1       & a_1   &   \cdots & 0    & 0         & 0\\
0        & c_2   & \lambda (x) - b_2    &   \cdots & 0    & 0         & 0 \\
&          \cdots    &     &         &      &           &   \\
0        & 0         & 0         &  \cdots &  c_{N-1} & \lambda (x) - b_{N-1} & a_{N-1}  \\
0        & 0         & 0         &  \cdots &  0   & c_N        & 
\lambda (x) - b_N  \\ \end{array}
       \right |\bdeq \det(\lambda(x)+{\cal G}_{N+1} ) , \label{ansatz}
\end{equation} 
where 
\begin{equation} b_n=-a_n-c_n \qquad \hbox{\rm for all} \quad n. \label{ansatz2}
\end{equation} 
The coefficients $a_n$ and $c_n$ for Hahn polynomials are
\begin{eqnarray*}
a_n & = & {(n+\alpha+\beta+1)(n+\alpha+1)(N-n)\over(2n+\alpha+\beta+1)
(2n+\alpha+\beta+2)} \\
c_n & = & {n(n+\alpha+\beta+N+1)(n+\beta)\over(2n+\alpha+\beta)
(2n+\alpha+\beta+1)}, 
\end{eqnarray*}
with $\lambda(x)=-x$, while the coefficients for Racah polynomials are
\begin{eqnarray*}
a_n & = & {(n+\alpha+1)(n+\alpha+\beta +1)(n+\gamma+1)(N-n)\over
(2n+\alpha+\beta+1)(2n+\alpha+\beta+2)} \\
c_n & = & {-n(n+\alpha+\beta+N+1)(n+\alpha+\beta-\gamma)(n+\beta)
\over(2n+\alpha+\beta)(2n+\alpha+\beta+1)},
\end{eqnarray*}
with $\beta+\gamma+1=-N$ and $\lambda(x)=-x(x+\gamma+\delta)$,
and the coefficients for $q$-Racah polynomials are
\begin{eqnarray}
a_n & = & {(1-abq^{n+1})(1-q^{n+1})(1-q^{n-N})(1-cq^{n+1})\over
(1-abq^{2n+1})(1-abq^{2n+2})}\label{as} \\
c_n & = & {cq^{-N}\over b} 
{(1-q^n)(1-bq^n)(1-abc^{-1}q^n)(1-abq^{n+N+1}) \over
(1-abq^{2n}) (1-abq^{2n+1})}, \label{cs} 
\end{eqnarray}
with $bdq=q^{-N}$ and $\lambda(x)=-(-q^{-x})(1-q^{x+1}cd)$ (see~\cite[Sec.~4]{A}). 

Recall that Hahn and Racah polynomials are just limiting cases of
$q$-Racah polynomials~\cite[Sec.~4]{A}. Precisely, if we denote Hahn polynomials 
by $H_{N+1}(\lambda(x);\alpha, \beta,N)$, Racah polynomials by
$RA_{N+1}(\lambda(x);\alpha,\beta,\gamma,N)$, and $q$-Racah polynomials
by $\qRa_{N+1}(\lambda(x);q,a,b,c,N)$, then we see that
\begin{eqnarray*}
RA_{N+1}(\lambda(x);\alpha,\beta,\gamma,N) & = & \lim_{q\to 1}
(-1)^{N+1}{\qRa_{N+1}(\lambda(x);q,q^\alpha,q^\beta,q^\gamma,N)\over
(1-q)^{2N+2}}, \\
H_{N+1}(\lambda(x);\alpha,\beta,N) & = & \lim_{\gamma\to \infty}
{RA_{N+1}(\lambda(x);\alpha,\beta,\gamma,N)\over\gamma^{N+1}}.
\end{eqnarray*}
So, it is enough to evaluate $q$-Racah polynomials.  The formula is given 
in~\cite{A} in the form
\[ \qRa_{N+1}(\lambda(x);q,a,b,c,N)=(-1)^{N+1}(q^{-x};q)_{N+1} (q^{x+1}cd;q)_{N+1}, \]
where $(\; ;\:)_k$ is the $q$-analogue of the shifted 
factorial~\cite[formula~(4.13)]{A}:
\[ (a;q)_k\eqbd \prod_{j=0}^{k-1} (1-aq^j). \]
Its equivalent form that is more suitable for an inductive proof is
\begin{equation}
 \qRa_{N+1}(\lambda(x);q,a,b,c,N)=\prod_{n=0}^N (\lambda(x)-\lambda(n)). 
\label{qRa} \end{equation}
This formula will be proved once we show that
\begin{equation}
 \qRa_{N+1}(\lambda(x);q,a,b,c,N)=\lambda(x)q^{-N}
\qRa_N(\widetilde\lambda(x-1); q,aq,b,cq,N-1), \label{qRaind} 
\end{equation}
where $\widetilde\lambda(x)$ corresponds to the parameters $(q,aq,b,cq,N-1)$
so that\begin{equation}
 \widetilde\lambda(x-1)=-(1-q^{-x+1}) (1-{cq\over b q^N}q^x)=
q\left(\lambda(x)+1+{c\over bq^N}-{1\over q}-{cq\over bq^N} \right). \label{connect}
 \end{equation}

To prove~(\ref{qRaind}), let us start with an observation about our ansatz
matrices ${\cal G}_{N+1}$ satisfying~(\ref{ansatz}) and~(\ref{ansatz2}). 
Suppose that such a ${\cal G}_{N+1}$ is transformed using
the matrix ${\cal T}_{N+1}$ given in~(\ref{t}). Then, as we already saw,  
${\cal T}_{N+1}^{-1} {\cal G}_{N+1} {\cal T}_{N+1}$ is block upper triangular: 
\[ \left[ \begin{array}{cc} 0 & * \\ 0 & {\cal M}_N
\end{array} \right]. \] 
Next, the transformation ${\cal S}_N$ given by~(\ref{s}) reduces ${\cal M}_N$
to the tridiagonal form
\begin{equation} {\cal S}_N {\cal M}_N  {\cal S}_N^{-1}= \left[    \begin{array}{cccccccc}
a_0+c_1        & a_1         & 0         & 0   &  \cdots & 0    & 0         & 0\\
c_1  & a_1+c_2       & a_2   & 0  &  \cdots & 0    & 0         & 0\\
0        & c_2   & a_2+c_3    & a_3   &  \cdots & 0    & 0         & 0 \\
&           & \cdots    &     &         &      &           &   \\
0        & 0         & 0         & 0   &  \cdots &  c_{N-2} & a_{N-2}+c_{N-1} & a_{N-1}  \\
0        & 0         & 0         & 0   &  \cdots &  0 & c_{N-1} & a_{N-1}+c_N  \\ 
\end{array}  \right ]. \label{tmp}   \end{equation}

Proving~(\ref{qRaind}) therefore amounts to showing that the matrix~(\ref{tmp}), 
with the entries $a_n$ and $c_n$ coming from the determinant 
$\qRa_{N+1}(\lambda(x);q,a,b,c,N)$, is similar to the matrix 
${1\over q}({\cal G}_N+\tilde\lambda(x-1)I_N)-\lambda(x)I_N$, where 
${\cal G}_N$ is an $N{\times}N$-matrix of the form~(\ref{ansatz})--(\ref{ansatz2}), 
with the entries $a_n$ and $c_n$ coming from the determinant 
$\qRa_N(\lambda(x-1);q,aq,b,cq,N-1)$. (Note that, due to~(\ref{connect}), the 
difference $\lambda(x)-{1\over q} \tilde\lambda(x-1)$ does not depend on $x$.)
The needed similarity   is realized by the diagonal matrix 
\[ \Lambda_N\eqbd \diag\left( (1-abq^2), {1\over q} (1-abq^4), {1\over q^2} (1-abq^6), 
\ldots,  {1\over q^{N-1}} (1-abq^{2N}) \right).  \] 
so that 
\begin{equation}
\Lambda_N^{-1} {\cal S}_N {\cal M}_N  {\cal S}_N^{-1}\Lambda_N={1\over q}({\cal G}_N+\tilde\lambda(x-1)I_N)-
\lambda(x)I_N.  \label{last}
\end{equation} 
Verification of~(\ref{last}) is straightforward for off-diagonal entries.
For diagonal entries, it reduces to verification of the identity
\[\left( a_n+c_{N+1}-1 -{c\over bq^N}+{1\over q}+{cq\over bq^N}\right) q=
a_{n+1}{1-abq^{2n+4}\over1-abq^{2n+2}}+q^2 c_n {1-abq^{2n}\over1-abq^{2n+2}}, 
\]
where $a_n$ are given by~(\ref{as}) and $c_n$ by~(\ref{cs}). This last identity can 
be checked using MATLAB Symbolic Math Toolbox. This finishes the proof of~(\ref{qRaind})
and~(\ref{qRa}).  \eproof

The determinant for Racah polynomials is therefore 
\[ RA_{N+1}(\lambda(x);\alpha,\beta,\gamma,N)=(-x)_{N+1}(x+\gamma+\delta+1)_{N+1} \]
and the determinant for Hahn polynomials is
\[H_{N+1}(\lambda(x);\alpha,\beta,N)=(-x)_{N+1}.  \]

\end{document}